\documentclass[12pt]{article}
\usepackage{amssymb,amsmath,amsthm,amsfonts,enumerate, bbm, mathdots}

\usepackage{latexsym}
\usepackage{amscd}
\usepackage{stmaryrd}
\usepackage{color}
\usepackage[all]{xypic}
\usepackage{epsfig}
\usepackage{graphics}
\usepackage{ifthen}
\usepackage{varioref}
\usepackage{rotating}
\usepackage{extarrows}
\usepackage{cite}
\usepackage{mathrsfs}

\numberwithin{equation}{section}

\theoremstyle{plain}
\newtheorem{thm}{Theorem}[section]
\newtheorem{cor}[thm]{Corollary}
\newtheorem{lem}[thm]{Lemma}

\newtheorem{rem}[thm]{Remark}

\definecolor{darkgreen}{rgb}{0.0625,0.64,0.0625}
\usepackage[
pdfstartview=FitH,
CJKbookmarks=true,
bookmarksnumbered=true,
bookmarksopen=true,
colorlinks,
pdfborder=001,
linkcolor=blue,
anchorcolor=green,
citecolor=red
]{hyperref}

\newfont{\scyr}{wncyr10 scaled 550}

\def\proof{\noindent {\bf Proof.\;}}

\def\sign{\operatorname{sign}}

\allowdisplaybreaks

\begin{document}
	
	\title{Generating functions of multiple $t$-star values of general level}
	
	\date{\small ~ \qquad\qquad School of Mathematical Sciences, Tongji University \newline No. 1239 Siping Road, Shanghai 200092, China}
	
	\author{Zhonghua Li\thanks{Email address: zhonghua\_li@tongji.edu.cn} ~and ~Lu Yan\thanks{Email address: 1910737@tongji.edu.cn}}

	\maketitle
	
	\begin{abstract}
		In this paper, we study the explicit expressions of multiple $t$-star values with an arbitrary number of blocks of twos of general level. We give an expression of a generating function of such values, which generalizes the results for multiple zeta-star values and multiple $t$-star values.  This derived generating function can provide expressions of multiple $t$-star values of general level in terms of the alternating multiple $t$-half values of general level with additional factorial and pochhammer symbol. As applications, some specific evaluations of multiple $t$-star values of general level with one-two-three or more general indices are given. These evaluations contribute to a deeper understanding of the properties of multiple $t$-star values of general level.
	\end{abstract}
	
	{\small
		{\bf Keywords} multiple $t$-star value, multiple $t$-star value of general level, alternating multiple $t$-half value, generating function.
	}
	
	{\small
		{\bf 2020 Mathematics Subject Classification} 11M32, 05A15.
	}
	
	
	\section{Introduction}\label{Sec:Intro}

In this paper, we study the explicit expressions of multiple $t$-star values with an arbitrary number of blocks of twos of general level. Let $N\in\mathbb{N}$ be fixed, where $\mathbb{N}$ is the set of positive integers. In \cite{Yuan-Zhao}, H. Yuan and J. Zhao introduced the multiple zeta values of level $N$. For any $\boldsymbol{s}=(s_1,\ldots,s_u)\in\mathbb{N}^u$ with $s_1>1$ and any $\boldsymbol{a}=(a_1,\ldots,a_u)\in {R_{N}}^u$ with $R_{N}=\mathbb{Z}/N\mathbb{Z}$, the multiple zeta value of level $N$ is defined by
		\begin{align*}
			\zeta_N(\boldsymbol{s};\boldsymbol{a})=\sum_{k_1>\cdots>k_u>0 \atop k_i\equiv a_i\mod N}
			\frac{1}{k_1^{s_1}\cdots k_u^{s_u}}.
		\end{align*}	
One can define the multiple zeta-star value of level $N$ by
		\begin{align*}
			\zeta_N^{\star}(\boldsymbol{s};\boldsymbol{a})=\sum_{k_1\geq\cdots\geq k_u>0 \atop k_i\equiv a_i\mod N}
			\frac{1}{k_1^{s_1}\cdots k_u^{s_u}}.
		\end{align*}	
		Note that if $a_1=\cdots=a
		_u$, we obtain the multiple $t$-(star) values of level $N$, which were studied by Z. Li and Z. Wang in \cite{Li-Wang}. For $\boldsymbol{s}=(s_1,\ldots,s_u)\in\mathbb{N}^u$ with $s_1>1$ and any $a\in \{1,2,\ldots,N\}$, the multiple $t$-value of level $N$ and the multiple $t$-star value of level $N$ are defined by
		\begin{align*}
			t_{N,a}(\boldsymbol{s})&=t_{N,a}(s_1,\ldots,s_u)=\sum_{k_1>\cdots>k_u\geq0}\frac{1}{{(Nk_1+a)}^{s_1}\cdots {(Nk_u+a)}^{s_u}},\\
			t^{\star}_{N,a}(\boldsymbol{s})&=t^{\star}_{N,a}(s_1,\ldots,s_u)=\sum_{k_1\geq\cdots\geq k_u\geq0}\frac{1}{{(Nk_1+a)}^{s_1}\cdots {(Nk_u+a)}^{s_u}},
		\end{align*}
		respectively. Here we treat the values $t_{N,a}(\emptyset)=t^{\star}_{N,a}(\emptyset)=1$. Setting $N=a=1$, we get the multiple zeta value and the multiple zeta-star value(cf. \cite{Hoffman1992,DZagier1994})
		 $$\zeta(\boldsymbol{s})=t_{1,1}(\boldsymbol{s})=\sum_{k_1>\cdots>k_u>0}
		 \frac{1}{k_1^{s_1}\cdots k_u^{s_u}},$$
		 $$\zeta^{\star}(\boldsymbol{s})=t_{1,1}^{\star}(\boldsymbol{s})=\sum_{k_1\geq\cdots\geq k_u>0}
		 \frac{1}{k_1^{s_1}\cdots k_u^{s_u}}.$$
		 Moreover,
		 $$t(\boldsymbol{s})=t_{2,1}(\boldsymbol{s})=\sum_{k_1>\cdots>k_u>0 \atop k_i:\text{odd}} \frac{1}{k_1^{s_1}\cdots k_u^{s_u}},$$ $$t^{\star}(\boldsymbol{s})=t_{2,1}^{\star}(\boldsymbol{s})=\sum_{k_1\geq\cdots\geq k_u>0 \atop k_i: \text{odd}} \frac{1}{k_1^{s_1}\cdots k_u^{s_u}},$$
		which are the multiple $t$-value and the multiple $t$-star value \cite{Hoffman2019} respectively. 

In \cite{DZagier2012}, D. Zagier provided explicit evaluation formulas for $\zeta(\{2\}^a,3,\{2\}^b)$ and $\zeta^{\star}(\{2\}^{a},3,\{2\}^{b})$ by establishing the generating functions, where $a,b\in\mathbb{N}_0=\mathbb{N}\cup \{0\}$. For an integer $c$, the $\{s\}^c$ denotes a sequence created by repeating $s$ exactly $c$ times if $c>0$; otherwise, $\{s\}^c$ is an empty string. Subsequently, additional proofs of Zagier's evaluation formulas have emerged; refer to, for instance \cite{PPT2014,Lee-Peng2018,Li2013,Li-Lupu-Orr}. Particularly in 2017, Kh. Hessami Pilehrood and T. Hessami Pilehrood \cite{PP2017} presented an alternative proof of Zagier's theorem, employing a new representation of the generating function of $\zeta^{\star}(\{2\}^{a},3,\{2\}^{b})$. Later in \cite{PP2019}, Kh. Hessami Pilehrood and T. Hessami Pilehrood studied the generating function of multiple zeta-star values with an arbitrary number of blocks of twos. Additionally, as outlined in \cite{PP2022}, they proposed an alternative approach to derive some established evaluations of multiple zeta-star values. Furthermore, they presented new explicit evaluations of $\zeta^{\star}(\{\{2\}^a,3,\{2\}^a,1\}^d)$ and $\zeta^{\star}(\{\{2\}^a,3,\{2\}^a,1\}^d,\{2\}^{a+1})$, with $a\in\mathbb{N}_0$ and $d\in\mathbb{N}$. It is noteworthly that there have been numerous contributions on evaluations of multiple zeta values and multiple zeta-star values, for example, see \cite{BBBL1998,BBBL2001, PPT2014, SMuneta2008, Ohno-Zudilin, DZagier1994, Zhao2016} and the references therein.
		
		On the other hand, employing the method of establishing the generating functions, T. Murakami\cite{TMurakami} obtained the evaluation of $t(\{2\}^a,3,\{2\}^b)$ and S. Charlton\cite{SCharlton} provided the evaluation of $t(\{2\}^a,1,\{2\}^b)$. In a recent paper\cite{Li-Lupu-Orr}, L. Lai, C. Lupu and D. Orr provided another proof of Murakami's result of $t(\{2\}^a,3,\{2\}^b)$. In \cite{Li-Wang}, Z. Li and Z. Wang found evaluations of $t^{\star}(\{2\}^a,3,\{2\}^b)$ and $t^{\star}(\{2\}^a,1,\{2\}^b)$.
		Consistent with the approach in a prior study \cite{PP2019}, we\cite{Li-Yan} delved into the generating functions of multiple $t$-star values with an arbitrary number of blocks of twos. Consequently, we derived analogous evaluation formulas as in \cite{Li-Wang} for $t^{\star}(\{2\}^a,3,\{2\}^b)$ and $t^{\star}(\{2\}^a,1,\{2\}^b)$, and gave additional evaluations of multiple t-star values showcasing relationships between multiple t-star values and alternating multiple t-half values.

		In this paper, we consider the generating functions of multiple $t$-star values of level $N$ with an arbitrary number of blocks of twos, extending previous findings on multiple zeta-star values and multiple $t$-star values (\cite{PP2019,Li-Yan}). In Section \ref{Sec: Generating functions of the finite form of multiple $t$-star values of level $N$}, we study the generating functions of the finite form of multiple $t$-star values of level $N$. Subsequently, in Section \ref{Sec: Multiple $t$-star values of level $N$ and their generating functions}, we provide explicit expressions of the generating functions of multiple $t$-star values of level $N$, and establish a formula of multiple $t$-star values of level $N$. Finally, Section \ref{Sec: Evaluations for MtSVs of level $N$} is dedicated to presenting various evaluations of multiple $t$-star values of level $N$.
	
	Throughout this paper, the following notations will be consistently employed. For any $\boldsymbol{s}=(s_1,\ldots,s_u)\in\mathbb{N}^u$ and $k, m \in\mathbb{N}_0$, we define that
	\begin{equation}\label{definitionofw}
		W_{k,m}^{\#}(\boldsymbol{s})=
		\begin{cases}
			\sum_{k\geq l_1\geq\cdots\geq l_u\geq m}\frac{2^{\triangle(k,l_1)+\triangle(l_1,l_2)+\cdots+\triangle(l_u,m)}}{(Nl_1+a)^{s_1}\cdots(Nl_u+a)^{s_u}} &\text{if\;} \boldsymbol{s}\neq\emptyset \text{\;and\;}k\geq m,\\
			&\\
			2^{\triangle(k,m)} &\text{otherwise},
		\end{cases}	
	\end{equation}
	where
	\[
	\triangle(k,m)=
	\begin{cases}
		0 &\text{if\;} k=m,\\
		1 &\text{if\;} k\neq m.
	\end{cases}
	\]
	For any $c\in\mathbb{N}_0$, define
	$$\delta(c)=
	\begin{cases}
		2 &\text{if\;} c=0,\\
		1 &\text{if\;} c=1,\\
		0 &\text{if\;} c\geq3.
	\end{cases}$$

	
	\section{The finite form of multiple $t$-star values of level $N$}\label{Sec: Generating functions of the finite form of multiple $t$-star values of level $N$}
	
	In this section, we prove a finite version of the generating function of multiple $t$-star values of level $N$. Let $n$ be a fixed nonnegative integer. For any $\boldsymbol{s}=(s_1,\ldots,s_u)\in\mathbb{N}^u$, we define the finite form of multiple $t$-star value of level $N$ by
	\begin{align*}
		_nt_{N,a}^{\star}(\boldsymbol{s})= \sum_{n\geq k_1\geq\cdots\geq k_u\geq0} \frac{1}{ (Nk_1+a)^{s_1}\cdots(Nk_u+a)^{s_u}},
	\end{align*}
	where $N$ is a fixed positive integer and $a\in\{1,2,\ldots,N\}$. Here we treat $_nt_{N,a}^{\star}(\emptyset)=1$. For $d\in\mathbb{N}_0$, $\boldsymbol{c}=(c_1,\ldots,c_d)\in\mathbb{N}^d$ and $\boldsymbol{z}=(z_0,z_1,\ldots,z_d)\in\mathbb{C}^{d+1}$, we denote the generating function of the finite form of multiple $t$-star values of level $N$ with an arbitrary number of blocks of twos as
	\begin{align*}
		D_n(\boldsymbol{c};\boldsymbol{z})=\sum_{b_0,b_1,\ldots,b_d\geq0} {}_nt_{N,a}^{\star}(\{2\}^{b_0},c_1,\{2\}^{b_1},\ldots,c_d,\{2\}^{b_d})z_0^{2b_0}z_1^{2b_1}\cdots z_d^{2b_d}.
	\end{align*}
	Recall that for an integer $c$, the $\{s\}^c$ denotes a sequence obtained by repeating $s$ exactly $c$ times if $c>0$; otherwise, $\{s\}^c$ is an empty string.
	
	 We give some lemmas in Subsection \ref{SubSec:PreLemThm2-1}, and prove an explicit expression of the generating function $D_n(\boldsymbol{c};\boldsymbol{z})$ in Subsection \ref{SubSec:ProfThm2-1}.

	\subsection{Preliminary Lemmas}\label{SubSec:PreLemThm2-1}

\begin{lem}
	For integers $a, N$ with $1\leq a \leq N$ and $n, l\in\mathbb{N}_0$, we have
	\begin{align}\label{l+1-n}
		\sum_{k=l+1}^n\frac{2\left(k+\frac{a}{N}\right)}{(n-k)!\left(\frac{2a}{N}\right)_{n+k+1}}=\frac{n-l}{(n-l)!\left(\frac{2a}{N}\right)_{n+l+1}}
	\end{align}
	and
	\begin{align}\label{alter-l+1-n}
		\sum_{k=l+1}^n\frac{2(-1)^k}{(n-k)!\left(\frac{2a}{N}\right)_{n+k+1}}=\frac{(-1)^{l+1}(n-l)}{\left(n+\frac{a}{N}\right)(n-l)!\left(\frac{2a}{N}\right)_{n+l+1}}.
	\end{align}
\end{lem}

\proof
We may assume that $n\geq l+1$. To prove \eqref{l+1-n},  we observe that
\begin{align}\label{I-recursion}
	\frac{2\left(k+\frac{a}{N}\right)}{(n-k)!\left(\frac{2a}{N}\right)_{n+k+1}}=I(n;k+1)-I(n;k)
\end{align}
for any integer $k$ with $l+1\leq k\leq n$, where
\begin{align*}
	I(n;j)=\frac{-1}{(n-k)!\left(\frac{2a}{N}\right)_{n+k}}
\end{align*}
for $l+1\leq j\leq n$	and $I(n;n+1)=0$. Then summing both sides of \eqref{I-recursion} over $k$ from $l+1$ to $n$, \eqref{l+1-n} follows easily.

Similarly, we set
\begin{align*}
	J(n;j)=\frac{(-1)^{j+1}}{\left(n+\frac{a}{N}\right)(n-j)!\left(\frac{2a}{N}\right)_{n+j}}
\end{align*}
for $l+1\leq j\leq n$	and $J(n;n+1)=0$. It is easy to see that for any integer $k$ with $l+1\leq k\leq n$,
\begin{align*}
	\frac{2(-1)^k}{(n-k)!\left(\frac{2a}{N}\right)_{n+k+1}}=J(n;k+1)-J(n;k).
\end{align*}
Then summing both sides of the above equation over $k$ from $l+1$ to $n$, we obtain \eqref{alter-l+1-n}.
\qed

\begin{lem}\label{alter-l-n-W}
	For integers $a, N$ with $1\leq a \leq N$ and $n,l, c\in\mathbb{N}_0$, we have
	\begin{align}\label{l-n-V}
		\sum_{k=l}^n\frac{(-1)^k}{(n-k)!\left(\frac{2a}{N}\right)_{n+k+1}}W_{k,l}^{\#}(\{1\}^c)=\frac{(-1)^l(Nl+a)}{(Nn+a)^{c+1}(n-l)!\left(\frac{2a}{N}\right)_{n+l+1}}.
	\end{align}
\end{lem}

\proof
We may assume that $n\geq l$. The proof is by induction on $c$. For $c=0$, we need to show
\begin{align}\label{c=0}
	\sum_{k=l}^n\frac{(-1)^k}{(n-k)!\left(\frac{2a}{N}\right)_{n+k+1}}\cdot2^{\triangle(k,l)}
	=\frac{(-1)^l(Nl+a)}{(Nn+a)(n-l)!\left(\frac{2a}{N}\right)_{n+l+1}}.
\end{align}
By the definition of $\triangle(k,l)$, the left-hand side of \eqref{c=0} is
$$\frac{(-1)^l}{(n-l)!\left(\frac{2a}{N}\right)_{n+l+1}}+\sum_{k=l+1}^n\frac{2(-1)^k}{(n-k)!\left(\frac{2a}{N}\right)_{n+k+1}}.$$
Applying \eqref{alter-l+1-n}, we find the left-hand side of \eqref{c=0} becomes
$$\frac{(-1)^l}{(n-l)!\left(\frac{2a}{N}\right)_{n+l+1}}+\frac{(-1)^{l+1}(n-l)}{\left(n+\frac{a}{N}\right)(n-l)!\left(\frac{2a}{N}\right)_{n+l+1}},$$
which exactly is the right-hand side of \eqref{c=0}.

If $c\geq1$, the left-hand side of \eqref{l-n-V} is
\begin{align*}
	&\sum_{k=l}^n\frac{(-1)^k}{(n-k)!\left(\frac{2a}{N}\right)_{n+k+1}}\sum_{k\geq l_1\geq\cdots\geq l_c\geq l}\frac{2^{\triangle(k,l_1)+\triangle(l_1,l_2)+\cdots+\triangle(l_c,l)}}{(Nl_1+a)\cdots(Nl_c+a)}\\
	=&\sum_{n\geq l_1\geq\cdots\geq l_c\geq l}\left\{\sum_{k=l_1}^{n}\frac{(-1)^k2^{\triangle(k,l_1)}}{(n-k)!\left(\frac{2a}{N}\right)_{n+k+1}}\right\}
	\frac{2^{\triangle(l_1,l_2)+\cdots+\triangle(l_c,l)}}{(Nl_1+a)(Nl_2+a)\cdots(Nl_c+a)}.
\end{align*}
Using \eqref{c=0} to deal with the part in the curly brace, we find the left-hand side of \eqref{l-n-V} is
$$\sum_{n\geq l_1\geq\cdots\geq l_c\geq l}\frac{(-1)^{l_1}(Nl_1+a)}{(Nn+a)(n-l_1)!\left(\frac{2a}{N}\right)_{n+l_1+1}}
\frac{2^{\triangle(l_1,l_2)+\cdots+\triangle(l_c,l)}}{(Nl_1+a)(Nl_2+a)\cdots(Nl_c+a)},$$
which is
$$\frac{1}{Nn+a}\sum_{l_1=l}^n\frac{(-1)^{l_1}}{(n-l_1)!\left(\frac{2a}{N}\right)_{n+l_1+1}}W_{l_1,l}^{\#}(\{1\}^{c-1}).$$
Then formula \eqref{l-n-V} easily follows from the inductive hypothesis for $c-1$.
\qed

	\subsection{Explicit expression of $D_n(\boldsymbol{c};\boldsymbol{z})$} \label{SubSec:ProfThm2-1}
	
		\begin{thm}\label{D_n-main-thm}
		Let $a, N$ be integers with $1\leq a \leq N$. For any $n, d\in\mathbb{N}_0$, $\boldsymbol{c}=(c_1,\ldots,c_d)\in(\mathbb{N}\setminus\{2\})^d$ and $\boldsymbol{z}=(z_0,z_1,\ldots,z_d)\in\mathbb{C}^{d+1}$ with $|z_j|<1, j=0,1,\ldots,d$, we have
		\begin{align}\label{D_n-main-formula}
			&D_n(\boldsymbol{c};\boldsymbol{z})=\sum_{n\geq k_0\geq k_1\geq\cdots\geq k_d\geq0}\frac{N\left[\left(\frac{a}{N}\right)_{n+1}\right]^2}{(n-k_0)!\left(\frac{2a}{N}\right)_{n+k_0+1}}\cdot\frac{\left(\frac{2a}{N}\right)_{k_d}}{k_d!(Nk_d+a)}\nonumber\\
			&\quad\quad\quad\quad\quad\quad\quad\quad\quad\quad\quad\quad\quad\quad\times\prod_{i=0}^{d}\frac{(-1)^{k_i\delta_i}(Nk_i+a)^{\delta_i-1}}{(Nk_i+a)^2-z_i^2}W_{k_{i-1},k_i}^{\#}(\{1\}^{c_i-3}),
		\end{align}
		where $k_{-1}=-1$, $\delta_i=\delta(c_i)+\delta(c_{i+1})$ with $c_0=1$ and $c_{d+1}=0$, the Pochhammer symbol $(a)_n$ is defined by
		$$(a)_n=\frac{\Gamma(a+n)}{\Gamma(a)}=
		\begin{cases}
			1 &\text{if\;} n=0,\\
			a(a+1)\cdots(a+n-1) &\text{if\;} n>0.
		\end{cases}$$
	\end{thm}
	
	Applying power series expansion $\frac{1}{(2k-1)^2-z^2}=\sum\limits_{b=0}^\infty\frac{z^{2b}}{(2k-1)^{2b+2}}$ and comparing the coefficients of $z_0^{2b_0}z_1^{2b_1}\cdots z_d^{2b_d}$ in \eqref{D_n-main-formula}, we get the following corollary.
	
	\begin{cor}\label{_nt_N,a^star-corollary}
		Let $a, N$ be integers with $1\leq a \leq N$. For any $n, d\in\mathbb{N}_0$, $c_1,\ldots,c_d\in\mathbb{N}\setminus\{2\}$ and    $b_0,b_1,\ldots,b_d\in\mathbb{N}_0$, we have
		\begin{align*}
			&_nt^{\star}_{N,a}(\{2\}^{b_0},c_1,\{2\}^{b_1},\ldots,c_d,\{2\}^{b_d})
			=\sum_{n\geq k_0\geq k_1\geq\cdots\geq k_d\geq0}\frac{N\left[\left(\frac{a}{N}\right)_{n+1}\right]^2}{(n-k_0)!\left(\frac{2a}{N}\right)_{n+k_0+1}}\cdot\frac{\left(\frac{2a}{N}\right)_{k_d}}{k_d!(Nk_d+a)}\nonumber\\
			&\quad\quad\quad\quad\quad\quad\quad\quad\quad\quad\quad\quad\quad\quad\quad\quad\quad\quad\quad\quad\times\prod_{i=0}^{d}\frac{(-1)^{k_i\delta_i}}{(Nk_i+a)^{2b_i+3-\delta_i}}W_{k_{i-1},k_i}^{\#}(\{1\}^{c_i-3}),
		\end{align*}
		where $k_{-1}=-1$, $\delta_i=\delta(c_i)+\delta(c_{i+1})$ with $c_0=1$ and $c_{d+1}=0$, and $(a)_n$ is the Pochhammer symbol.
	\end{cor}
	
	\noindent \textbf{Proof of Theorem \ref{D_n-main-thm}.}
	We apply induction on $n+d$. For $n=0$, we have
	\begin{align*}
		D_0(\boldsymbol{c}; \boldsymbol{z})=\sum_{b_0,b_1,\ldots,b_d\geq0}\frac{z_0^{2b_0}\cdots z_d^{2b_d}}{a^{2b_0+\cdots+2b_d+c_1+\cdots+c_d}}=\frac{1}{a^{c_1+\cdots+c_d}}\prod_{i=0}^{d}\frac{a^2}{a^2-z_i^2},
	\end{align*}
	and the right-hand side of \eqref{D_n-main-formula} is
	\begin{align*}
		&\frac{1}{2}\prod_{i=0}^d\frac{a^{\delta_i-1}}{a^2-z_i^2}\cdot W_{-1,0}^{\#}(\emptyset)W_{0,0}^{\#}(\{1\}^{c_1-3})\cdots W_{0,0}^{\#}(\{1\}^{c_d-3})\\
		=&\prod_{i=0}^d\frac{a^2}{a^2-z_i^2}\cdot\frac{a^{2\delta(c_1)+\cdots+2\delta(c_d)-3d}}{a^{c_1-3+2\delta(c_1)}\cdots a^{c_d-3+2\delta(c_d)}}\\
		=&\frac{1}{a^{c_1+\cdots+c_d}}\prod_{i=0}^{d}\frac{a^2}{a^2-z_i^2}.
	\end{align*}
	So \eqref{D_n-main-formula} is true for $n=0$.
	
	For $d=0$, we observe that for $n>0$,
	$$_nt_{N,a}^{\star}(\{2\}^{b_0})=\sum_{k=0}^{b_0}\frac{_{n-1}t_{N,a}^{\star}(\{2\}^k)}{(Nn+a)^{2(b_0-k)}}.$$
	Hence we obtain a recursive formula
	$$D_n( ; z_0)=\frac{(Nn+a)^2}{(Nn+a)^2-z_0^2}D_{n-1}( ; z_0).$$
	Using the above formula repeatedly, we have
	\begin{align*}
		D_n( ; z_0)=\prod_{j=1}^{n}\frac{(Nj+a)^2}{(Nj+a)^2-z_0^2}D_{0}( ; z_0)
		=\prod_{j=0}^{n}\frac{(Nj+a)^2}{(Nj+a)^2-z_0^2}.
	\end{align*}
	By the method of partial fractional decomposition, we get
	\begin{align*}
		D_n( ; z_0)
		=&\sum_{k=0}^{n}\frac{2(-1)^{k}N\left[\left(\frac{a}{N}\right)_{n+1}\right]^2\left(\frac{2a}{N}\right)_{k}}{k!(n-k)!\left(\frac{2a}{N}\right)_{n+k+1}(Nk+a)}\cdot\frac{(Nk+a)^2}{(Nk+a)^2-z_0^2}\\
		=&\sum_{k=0}^{n}\frac{N\left[\left(\frac{a}{N}\right)_{n+1}\right]^2}{(n-k)!\left(\frac{2a}{N}\right)_{n+k+1}}\cdot\frac{\left(\frac{2a}{N}\right)_{k}}{k!(Nk+a)}\cdot\frac{2(-1)^k(Nk+a)^2}{(Nk+a)^2-z_0^2}.
	\end{align*}
	Therefore, \eqref{D_n-main-formula} is proved for $d=0$.	
	
	Now assume that $n>0$ and $d>0$. Since
	\begin{align*}
		_nt^{\star}_{N,a}(\{2\}^{b_0},c_1,\{2\}^{b_1},\ldots,c_d,\{2\}^{b_d})
		=&\sum_{k=0}^{b_0}\frac{1}{(Nn+a)^{2b_0-2k}}{}_{n-1}t^{\star}_{N,a}(\{2\}^{k},c_1,\{2\}^{b_1},\ldots,c_d,\{2\}^{b_d})\\
		&+\frac{1}{(Nn+a)^{2b_0+c_1}}{}_nt^{\star}_{N,a}(\{2\}^{b_1},c_2,\{2\}^{b_2},\ldots,c_d,\{2\}^{b_d}),
	\end{align*}
	we deduce that
	\begin{align*}
		D_n(\boldsymbol{c}; \boldsymbol{z})
		=&\sum_{b_0,\ldots,b_d\geq0}\sum_{k=0}^{b_0}\frac{1}{(Nn+a)^{2b_0-2k}}	{}_{n-1}t^{\star}_{N,a}(\{2\}^{k},c_1,\{2\}^{b_1},\ldots,c_d,\{2\}^{b_d})z_0^{2b_0}\cdots z_d^{2b_d},\\
		&+\sum_{b_0\geq0}\frac{z_0^{2b_0}}{(Nn+a)^{2b_0+c_1}}D_n(\boldsymbol{c}^-; \boldsymbol{z}^-),
	\end{align*}
	where $\boldsymbol{c}^-=(c_2,\ldots,c_d)$ and $\boldsymbol{z}^-=(z_1,z_2,\ldots,z_d)$. Then it is easy to get the following recursive formula
	\begin{align}\label{recurence-D_n}
		D_n(\boldsymbol{c}; \boldsymbol{z})=\frac{(Nn+a)^2}{(Nn+a)^2-z_0^2}D_{n-1}(\boldsymbol{c}; \boldsymbol{z})+\frac{(Nn+a)^{2-c_1}}{(Nn+a)^2-z_0^2}D_n(\boldsymbol{c}^-; \boldsymbol{z}^-).
	\end{align}
	Let $\widetilde{D}_n(\boldsymbol{c};\boldsymbol{z})$ denote the right-hand side of \eqref{D_n-main-formula}. According to \eqref{recurence-D_n}, we need to prove
	\begin{align}\label{D_n-main-new}
		\frac{(Nn+a)^2}{(Nn+a)^2-z_0^2}D_{n-1}(\boldsymbol{c}; \boldsymbol{z})-\widetilde{D}_n(\boldsymbol{c};\boldsymbol{z})
		=-\frac{(Nn+a)^{2-c_1}}{(Nn+a)^2-z_0^2}D_n(\boldsymbol{c}^-; \boldsymbol{z}^-).
	\end{align}
	By the induction hypothesis for $D_{n-1}(\boldsymbol{c}; \boldsymbol{z})$, the left hand-side of \eqref{D_n-main-new} is
	\begin{align*}
         &\frac{(Nn+a)^2}{(Nn+a)^2-z_0^2}\sum_{n-1\geq k_0\geq k_1\geq\cdots\geq k_d\geq0}\frac{N\left[\left(\frac{a}{N}\right)_{n}\right]^2}{(n-1-k_0)!\left(\frac{2a}{N}\right)_{n+k_0}}\cdot\frac{\left(\frac{2a}{N}\right)_{k_d}}{k_d!(Nk_d+a)}\\
         &\quad\quad\quad\quad\quad\quad\quad\quad\quad\quad\quad\times\prod_{i=0}^{d}\frac{(-1)^{k_i\delta_i}(Nk_i+a)^{\delta_i-1}}{(Nk_i+a)^2-z_i^2}W_{k_{i-1},k_i}^{\#}(\{1\}^{c_i-3})\\
         &-\sum_{n\geq k_0\geq k_1\geq\cdots\geq k_d\geq0}\frac{N\left[\left(\frac{a}{N}\right)_{n+1}\right]^2}{(n-k_0)!\left(\frac{2a}{N}\right)_{n+k_0+1}}\cdot\frac{\left(\frac{2a}{N}\right)_{k_d}}{k_d!(Nk_d+a)}\nonumber\\
         &\quad\quad\quad\quad\quad\quad\quad\quad\quad\quad\quad\quad\quad\quad\times\prod_{i=0}^{d}\frac{(-1)^{k_i\delta_i}(Nk_i+a)^{\delta_i-1}}{(Nk_i+a)^2-z_i^2}W_{k_{i-1},k_i}^{\#}(\{1\}^{c_i-3}),
	\end{align*}
	which is
	\begin{align*}
		&\sum_{n\geq k_0\geq\cdots\geq k_d\geq0}\prod_{i=0}^{d}\frac{(-1)^{k_i\delta_i}(Nk_i+a)^{\delta_i-1}}{(Nk_i+a)^2-z_i^2}W_{k_{i-1},k_i}^{\#}(\{1\}^{c_i-3})\cdot\frac{N\left(\frac{2a}{N}\right)_{k_d}}{k_d!(Nk_d+a)}\\
		&\quad\quad\quad\times\left\{\frac{(Nn+a)^2}{(Nn+a)^2-z_0^2}\cdot\frac{\left[\left(\frac{a}{N}\right)_{n}\right]^2}{(n-1-k_0)!\left(\frac{2a}{N}\right)_{n+k_0}}
		-\frac{\left[\left(\frac{a}{N}\right)_{n+1}\right]^2}{(n-k_0)!\left(\frac{2a}{N}\right)_{n+k_0+1}}\right\}.
	\end{align*}	
	As the part in braces can be simplified as
	$$-\frac{\left[\left(\frac{a}{N}\right)_{n+1}\right]^2}{(n-k_0)!\left(\frac{2a}{N}\right)_{n+k_0+1}}\cdot\frac{(Nk_0+a)^2-z_0^2}{(Nn+a)^2-z_0^2},$$
	the left-hand side of \eqref{D_n-main-new} becomes
	\begin{align}\label{difference}
		&\frac{-2N\left[\left(\frac{a}{N}\right)_{n+1}\right]^2}{(Nn+a)^2-z_0^2}\sum_{n\geq k_0\geq k_1\geq\cdots\geq k_d\geq0}\frac{(-1)^{k_0(1+\delta(c_1))}(Nk_0+a)^{\delta(c_1)}}{(n-k_0)!\left(\frac{2a}{N}\right)_{n+k_0+1}}\cdot\frac{\left(\frac{2a}{N}\right)_{k_d}}{k_d!(Nk_d+a)}\nonumber\\
		&\quad\quad\quad\quad\quad\quad\quad\times\prod_{i=1}^{d}\frac{(-1)^{k_i\delta_i}(Nk_i+a)^{\delta_i-1}}{(Nk_i+a)^2-z_i^2}W_{k_{i-1},k_i}^{\#}(\{1\}^{c_i-3}).
	\end{align}
	Let $\sum\nolimits_{k_0}$ denote the inner sum over $k_0$ in \eqref{difference}, that is
	\begin{align*}
		\sum\nolimits_{k_0}=\sum_{k_0=k_1}^{n}\frac{(-1)^{k_0(1+\delta(c_1))}(Nk_0+a)^{\delta(c_1)}}{(n-k_0)!\left(\frac{2a}{N}\right)_{n+k_0+1}}W_{k_0,k_1}^{\#}(\{1\}^{c_1-3}).
	\end{align*}
	
	If $c_1=1$, $\delta(c_1)=1$,
	\begin{align*}
		\sum\nolimits_{k_0}=&\sum_{k_0=k_1}^{n}\frac{Nk_0+a}{(n-k_0)!\left(\frac{2a}{N}\right)_{n+k_0+1}}2^{\triangle(k_0,k_1)}.\\
		=&\frac{Nk_1+a}{(n-k_1)!\left(\frac{2a}{N}\right)_{n+k_1+1}}+2\sum_{k_0=k_1+1}^{n}\frac{Nk_0+a}{(n-k_0)!\left(\frac{2a}{N}\right)_{n+k_0+1}}.
	\end{align*}
	Then by \eqref{l+1-n}, we have
	$$\sum\nolimits_{k_0}=\frac{Nk_1+a}{(n-k_1)!\left(\frac{2a}{N}\right)_{n+k_1+1}}+\frac{Nn-Nk_1}{(n-k_1)!\left(\frac{2a}{N}\right)_{n+k_1+1}}=\frac{Nn+a}{(n-k_1)!\left(\frac{2a}{N}\right)_{n+k_1+1}}.$$
	Using the above result together with \eqref{difference}, we find the left-hand side of \eqref{D_n-main-new} equals
	\begin{align*}
		&\frac{-(Nn+a)}{(Nn+a)^2-z_0^2}\sum_{n\geq k_1\geq\cdots\geq k_d\geq0}\frac{N\left[\left(\frac{a}{N}\right)_{n+1}\right]^2
			}{(n-k_1)!\left(\frac{2a}{N}\right)_{n+k_1+1}}\cdot\frac{\left(\frac{2a}{N}\right)_{k_d}}{k_d!(Nk_d+a)}\\
		&\quad\quad\quad\quad\quad\quad\quad\times\frac{2(-1)^{k_1(1+\delta(c_2))}(Nk_1+a)^{\delta(c_2)}}{(Nk_1+a)^2-z_1^2}\prod_{i=2}^{d}\frac{(-1)^{k_i\delta_i}(Nk_i+a)^{\delta_i-1}}{(Nk_i+a)^2-z_i^2}W_{k_{i-1},k_i}^{\#}(\{1\}^{c_i-3}),
	\end{align*}
	which by the inductive hypothesis for $D_{n}(\boldsymbol{c}^-;\boldsymbol{z}^-)$ is
	$$-\frac{Nn+a}{(Nn+a)^2-z_0^2}D_{n}(\boldsymbol{c}^-;\boldsymbol{z}^-).$$
	Therefore, the theorem is proved in this case.
	
	If $c_1\geq3$, $\delta(c_1)=0$, by Lemma \ref{alter-l-n-W}, we have
	\begin{align*}
        \sum\nolimits_{k_0}=\sum_{k_0=k_1}^{n}\frac{(-1)^{k_0}}{(n-k_0)!\left(\frac{2a}{N}\right)_{n+k_0+1}}W_{k_0,k_1}^{\#}(\{1\}^{c_1-3})=\frac{(-1)^{k_1}(Nn+a)^{2-c_1}(Nk_1+a)}{(n-k_1)!\left(\frac{2a}{N}\right)_{n+k_1+1}}
	\end{align*}
	Combining this formula with \eqref{difference}, we obtain the left-hand side of \eqref{D_n-main-new} is
	\begin{align*}
		&\frac{-(Nn+a)^{2-c_1}}{(Nn+a)^2-z_0^2}\sum_{n\geq k_1\geq\cdots\geq k_d\geq0}\frac{N\left[\left(\frac{a}{N}\right)_{n+1}\right]^2
		}{(n-k_1)!\left(\frac{2a}{N}\right)_{n+k_1+1}}\cdot\frac{\left(\frac{2a}{N}\right)_{k_d}}{k_d!(Nk_d+a)}\\
		&\quad\quad\quad\quad\quad\quad\quad\times\frac{2(-1)^{k_1(1+\delta(c_2))}(Nk_1+a)^{\delta(c_2)}}{(Nk_1+a)^2-z_1^2}\prod_{i=2}^{d}\frac{(-1)^{k_i\delta_i}(Nk_i+a)^{\delta_i-1}}{(Nk_i+a)^2-z_i^2}W_{k_{i-1},k_i}^{\#}(\{1\}^{c_i-3}),
	\end{align*}
	which by the inductive hypothesis for $D_{n}(\boldsymbol{c}^-;\boldsymbol{z}^-)$ is
	$$-\frac{(Nn+a)^{2-c_1}}{(Nn+a)^2-z_0^2}D_{n}(\boldsymbol{c}^-;\boldsymbol{z}^-).$$
	Hence we conclude that $D_n(\boldsymbol{c}; \boldsymbol{z})=\widetilde{D}_n(\boldsymbol{c};\boldsymbol{z})$. The proof is completed.
	\qed


	\section{Multiple $t$-star values of level $N$ and their generating functions}\label{Sec: Multiple $t$-star values of level $N$ and their generating functions}
	
	In this section, we prove  the correctness of limit transfer from the generating functions of the finite form of multiple $t$-star values of level $N$ to the generating functions of multiple $t$-star values of level $N$.
	For $d\in\mathbb{N}_0$, $\boldsymbol{c}=(c_1,\ldots,c_d)\in{\mathbb{N}}^d$ with $c_1>1$ and $\boldsymbol{z}=(z_0,z_1,\ldots,z_d)\in\mathbb{C}^{d+1}$, we define the generating function
	\begin{align*}
		D(\boldsymbol{c}; \boldsymbol{z})=\sum_{b_0,b_1,\ldots,b_d\geq0}	t_{N,a}^{\star}(\{2\}^{b_0},c_1,\{2\}^{b_1},\ldots,c_d,\{2\}^{b_d})z_0^{2b_0}z_1^{2b_1}\cdots z_d^{2b_d}.
	\end{align*}
	
	In Subsection \ref{SubSec:LemsThms3-1,3,4}, we prepare some lemmas that will be used in the proof of Theorems \ref{D-main-thm} and \ref{t^star-levelN-theorem}. Then we give Theorems \ref{D-main-thm} and \ref{t^star-levelN-theorem} and prove them in Subsection \ref{SubSec:ProfThms3-4,5}.
	
	\subsection{Lemmas}\label{SubSec:LemsThms3-1,3,4}
	
	\begin{lem}\label{limit-D_n=D-lemma}
		Let $d\in\mathbb{N}_0$, $\boldsymbol{c}=(c_1,\ldots,c_d)\in(\mathbb{N}\setminus\{2\})^d$ with $c_1\geq3$, and let $\boldsymbol{z}=(z_0,z_1,\ldots,z_d)\in\mathbb{C}^{d+1}$
		with $|z_j|<1, j=0,1,\ldots,d$. Then
		\begin{align}\label{limit-D_n=D}
			\lim_{n\rightarrow\infty}D_n(\boldsymbol{c}; \boldsymbol{z})=D(\boldsymbol{c}; \boldsymbol{z}),
		\end{align}
		and the convergence is uniform in any closed region $E$: $|z_0|\leq u_0<1, |z_1|\leq u_1<1, \ldots, |z_d|\leq u_d<1$.
	\end{lem}
	
	\proof
	For $|z|<1$, we deduce that
	\begin{align}\label{sum-_nmt^*_Na}
		\prod_{k=m}^{n}\left(1-\frac{z^2}{(Nk+a)^2}\right)^{-1}=\sum_{l=0}^{\infty}{}_{n,m}t_{N,a}^{\star}(\{2\}^l)z^{2l}.
	\end{align}
	Here we define that, for integers $n\geq m >0$ and $\boldsymbol{s}=(s_1,\ldots, s_u)\in\mathbb{N}^u$,
	\begin{align*}
		_{n,m}t_{N,a}^{\star}(\boldsymbol{s})=_{n,m}t_{N,a}^{\star}(s_1,\ldots, s_u)=\sum_{n\geq k_1\geq\cdots\geq k_u\geq m} \frac{1}{ (Nk_1+a)^{s_1}\cdots(Nk_u+a)^{s_u}}.
	\end{align*}
	We also set $_{n,m}t^{\star}_{N,a}(\emptyset)=1$. If $s_1>1$, we let $n$ tend to infinity and define
	$$_{\infty,m}t_{N,a}^{\star}(\boldsymbol{s})=\sum_{k_1\geq\cdots\geq k_u\geq m} \frac{1}{ (Nk_1+a)^{s_1}\cdots(Nk_u+a)^{s_u}}.$$
	
	If $d=0$, using \eqref{sum-_nmt^*_Na}, we get
	\begin{align*}
		D_n( ;z_0)=\prod_{k=0}^n\left(1-\frac{z_0^2}{(Nk+a)^2}\right)^{-1}
\quad \text{and}\quad\quad
		D( ;z_0)=\prod_{k=0}^\infty\left(1-\frac{z_0^2}{(Nk+a)^2}\right)^{-1}.
	\end{align*}
	Hence the lemma is true for $d=0$.
	
	Now assume that $d\geq1$. Using \eqref{sum-_nmt^*_Na}, we have
	\begin{align*}
		&D_n(\boldsymbol{c}; \boldsymbol{z})\\
		=&\sum_{n\geq k_1\geq\cdots\geq k_d\geq0}\frac{\prod\limits_{k=k_1}^n\left(1-\frac{z_0^2}{(Nk+a)^2}\right)^{-1}\prod\limits_{k=k_2}^{k_1}\left(1-\frac{z_1^2}{(Nk+a)^2}\right)^{-1}\cdots\prod\limits_{k=0}^{k_d}\left(1-\frac{z_d^2}{(Nk+a)^2}\right)^{-1}}{(Nk_1+a)^{c_1}\cdots(Nk_d+a)^{c_d}}.
	\end{align*}
	Setting
	\begin{align*}
		&D_{n}^{\infty}(\boldsymbol{c}; \boldsymbol{z})\\
    	=&\sum_{n\geq k_1\geq\cdots\geq k_d\geq0}\frac{\prod\limits_{k=k_1}^{\infty}\left(1-\frac{z_0^2}{(Nk+a)^2}\right)^{-1}\prod\limits_{k=k_2}^{k_1}\left(1-\frac{z_1^2}{(Nk+a)^2}\right)^{-1}\cdots\prod\limits_{k=0}^{k_d}\left(1-\frac{z_d^2}{(Nk+a)^2}\right)^{-1}}{(Nk_1+a)^{c_1}\cdots(Nk_d+a)^{c_d}},
	\end{align*}
	then
	\begin{align}\label{Dn^inf-Dn}
		&|D_{n}^{\infty}(\boldsymbol{c}; \boldsymbol{z})-D_{n}(\boldsymbol{c}; \boldsymbol{z})|\nonumber\\
		\leq&\sum_{n\geq k_1\geq\cdots\geq k_d\geq0}\frac{\prod\limits_{k=k_1}^n\left(1-\frac{z_0^2}{(Nk+a)^2}\right)^{-1}\prod\limits_{k=k_2}^{k_1}\left(1-\frac{z_1^2}{(Nk+a)^2}\right)^{-1}\cdots\prod\limits_{k=0}^{k_d}\left(1-\frac{z_d^2}{(Nk+a)^2}\right)^{-1}}{(Nk_1+a)^{c_1}\cdots(Nk_d+a)^{c_d}}\nonumber\\
		&\quad\quad\quad\quad\quad\quad\quad\quad\quad\quad\quad\quad\quad\times\left|\prod\limits_{k=n+1}^{\infty}\left(1-\frac{z_0^2}{(Nk+a)^2}\right)^{-1}-1\right|.
	\end{align}
	Assume that $|z_0|\leq u_0<1, |z_1|\leq u_1<1, \ldots, |z_d|\leq u_d<1$. Using the infinite product formula for the sine function
	$$\frac{\sin\pi z}{\pi z}=\prod\limits_{k=1}^\infty\left(1-\frac{z^2}{k^2}\right)<\prod\limits_{k=0}^\infty\left(1-\frac{z^2}{(Nk+a)^2}\right),$$
	we obtain
	\begin{align}\label{Dn-upper-bound}
		&\sum_{n\geq k_1\geq\cdots\geq k_d\geq0}\frac{\prod\limits_{k=k_1}^n\left(1-\frac{z_0^2}{(Nk+a)^2}\right)^{-1}\prod\limits_{k=k_2}^{k_1}\left(1-\frac{z_1^2}{(Nk+a)^2}\right)^{-1}\cdots\prod\limits_{k=0}^{k_d}\left(1-\frac{z_d^2}{(Nk+a)^2}\right)^{-1}}{(Nk_1+a)^{c_1}\cdots(Nk_d+a)^{c_d}}\nonumber\\
		<&\prod_{i=0}^d\frac{\pi|z_i|}{\sin(\pi |z_i|)}{}_nt_{N,a}^{\star}(\boldsymbol{c})
		<\prod_{i=0}^d\frac{\pi u_i}{\sin(\pi u_i)}t_{N,a}^{\star}(\boldsymbol{c}).
	\end{align}
	Using \eqref{sum-_nmt^*_Na}, we have
	\begin{align*}
		\left|\prod\limits_{k=n+1}^{\infty}\left(1-\frac{z_0^2}{(Nk+a)^2}\right)^{-1}-1\right|=\left|\sum_{l=1}^{\infty}{}_{\infty,n+1}t_{N,a}^{\star}(\{2\}^l)z_0^{2l}\right|<\sum_{l=1}^{\infty}{}_{\infty,n+1}t_{N,a}^{\star}(\{2\}^l).
	\end{align*}
	Note that
	\begin{align*}
		_{\infty,n+1}t_{N,a}^{\star}(\{2\}^l)&=\sum_{k_1\geq\cdots\geq k_l\geq n+1}\frac{1}{(Nk_1+a)^2\cdots (Nk_l+a)^2}\\
		&<\left(\sum_{k=n+1}^{\infty}\frac{1}{(Nk+a)^2}\right)^l<\left(\int_{n}^{\infty}\frac{\mathrm{d}x}{(Nx+a)^2}\right)^l=\left(\frac{1}{N(Nn+a)}\right)^l,
	\end{align*}
	then we obtain
	\begin{align}\label{1/4n-3}
		\left|\prod\limits_{k=n+1}^{\infty}\left(1-\frac{z_0^2}{(Nk+a)^2}\right)^{-1}-1\right|<\sum_{l=1}^{\infty}\left(\frac{1}{N(Nn+a)}\right)^l=\frac{1}{N(Nn+a)-1}.
	\end{align}
	Using \eqref{Dn^inf-Dn}, \eqref{Dn-upper-bound} and \eqref{1/4n-3}, we have
	\begin{align*}
		\left|D_{n}^{\infty}(\boldsymbol{c}; \boldsymbol{z})-D_{n}(\boldsymbol{c}; \boldsymbol{z})\right|
		<\frac{1}{N(Nn+a)-1}\prod_{i=0}^d\frac{\pi u_i}{\sin(\pi u_i)}t_{N,a}^{\star}(\boldsymbol{c})\rightarrow0,
	\end{align*}
	as $n\rightarrow\infty$ on the closed region $E$. Similarly, we find
	\begin{align*}
		&\left|D(\boldsymbol{c}; \boldsymbol{z})-D_{n}^{\infty}(\boldsymbol{c}; \boldsymbol{z})\right|\\
		\leq&\sum_{k_1\geq\cdots\geq k_d\geq0\atop k_1>n}\frac{\prod\limits_{k=k_1}^{\infty}\left(1-\frac{|z_0|^2}{(Nk+a)^2}\right)^{-1}\prod\limits_{k=k_2}^{k_1}\left(1-\frac{|z_1|^2}{(Nk+a)^2}\right)^{-1}\prod\limits_{k=0}^{k_d}\left(1-\frac{|z_d|^2}{(Nk+a)^2}\right)^{-1}}{(Nk_1+a)^{c_1}\cdots(Nk_d+a)^{c_d}}\nonumber\\
		<&\prod_{i=0}^d\frac{\pi u_i}{\sin(\pi u_i)}\left(t_{N,a}^{\star}(\boldsymbol{c})-{}_nt_{N,a}^{\star}(\boldsymbol{c})\right)\rightarrow0,
	\end{align*}
	as $n\rightarrow\infty$ on the closed region $E$. Therefore, the proof is finished.
	\qed

	\begin{lem}\label{limit-M_k-lemma}
		Let $N\in\mathbb{N}$, $a\in\{1,2,\ldots,N\}$, $\varepsilon, p, q, C \in\mathbb{R}$ with $\varepsilon>0, p> 1, C>0$, and let $M_k$ be real numbers satisfying
		$$|M_k|< \frac{C\log^q(Nk+N+a)}{k^p+\varepsilon}$$
		for any $k\in\mathbb{N}_0$. Then we have
		\begin{align*}
			\lim_{n\rightarrow\infty}\sum_{k=0}^{n}|M_k|\frac{N\left[\left(\frac{a}{N}\right)_{n+1}\right]^2}{(n-k)!\left(\frac{2a}{N}\right)_{n+k+1}}=\sum_{k=0}^{\infty}|M_k|\frac{N}{B\left(\frac{a}{N},\frac{a}{N}\right)},
		\end{align*}
	where $B(a,b)$ is the Beta function.
	\end{lem}
	
	\proof
	Let $$U_{k,n}=\frac{N\left[\left(\frac{a}{N}\right)_{n+1}\right]^2}{(n-k)!\left(\frac{2a}{N}\right)_{n+k+1}} =\frac{N\Gamma\left(\frac{2a}{N}\right)}{\left[\Gamma\left(\frac{a}{N}\right)\right]^2}\cdot\frac{\left[\Gamma\left(n+1+\frac{a}{N}\right)\right]^2}{\Gamma(n+1-k)\Gamma\left(n+1+k+\frac{2a}{N}\right)},$$
	where $\Gamma(a)$ is the Gamma function. For any $k\in\mathbb{N}_0$, applying the well-known Stirling's formula (cf. \cite{Diaconis-Freedman,Feller})
	$$\Gamma(x)\sim\sqrt{2\pi}(x-1)^{x-\frac{1}{2}}e^{-(x-1)},\quad x\rightarrow\infty,$$
	and the limit formula $$\left(1+\frac{1}{x}\right)^x\sim e,\quad x\rightarrow\infty,$$
	we obtain that
	\begin{align*}
		&\frac{\left[\Gamma\left(n+1+\frac{a}{N}\right)\right]^2}{\Gamma(n+1-k)\Gamma\left(n+1+k+\frac{2a}{N}\right)}\\
		\sim&\frac{n+\frac{a}{N}}{\sqrt{(n-k)\left(n+k+\frac{2a}{N}\right)}}\left(1+\frac{k+\frac{a}{N}}{n-k}\right)^{n-k}\left(1+\frac{-k-\frac{a}{N}}{n+k+\frac{2a}{N}}\right)^{n+k+\frac{2a}{N}}\\
		\sim&1,\quad n\rightarrow\infty.
	\end{align*}
	Therefore, we have
	$$\lim_{n\rightarrow\infty}U_{k,n}=\frac{N\Gamma\left(\frac{2a}{N}\right)}{\left[\Gamma\left(\frac{a}{N}\right)\right]^2}=\frac{N}{B\left(\frac{a}{N},\frac{a}{N}\right)},$$	
	which leads to the conclusion that $U_{k,n}$ is a bounded sequence of numbers. It is easy to verify that $U_{k,n}$ is decreasing with respect to $k$, so we have
	$0\leq U_{k,n}\leq U_{0,n}.$
	Note that
	$$\frac{U_{k,n+1}}{U_{k,n}}=\frac{\left(n+1+\frac{a}{N}\right)^2}{(n-k+1)\left(n+k+1+\frac{2a}{N}\right)}=\frac{n^2+\left(2+\frac{2a}{N}\right)n+\left(1+\frac{a}{N}\right)^2}{n^2+\left(2+\frac{2a}{N}\right)n+(1-k)\left(k+1+\frac{2a}{N}\right)}>1,$$
	so $U_{k,n}$ is increasing with respect to $n$.
	Since $|M_k|< \frac{C\log^q(Nk+N+a)}{k^p+\varepsilon}$ with $p>1$,
	\begin{align}\label{n-infty1}
		\sum_{k=0}^{\infty}|M_k|U_{k,n}-\sum_{k=0}^{n}|M_k|U_{k,n}=\sum_{k=n+1}^{\infty}|M_k|U_{k,n}\rightarrow0
	\end{align}
	as $n\rightarrow\infty$. By Monotone Convergence Theorem, we have
	\begin{align}\label{n-infty2}
		\lim\limits_{n\rightarrow\infty}\sum_{k=0}^{\infty}|M_k|U_{k,n}=\sum_{k=0}^{\infty}|M_k|\lim\limits_{n\rightarrow\infty}U_{k,n}=\sum_{k=0}^{\infty}|M_k|\frac{N}{B\left(\frac{a}{N},\frac{a}{N}\right)}.
	\end{align}
	Hence, \eqref{n-infty1} and \eqref{n-infty2} can deduce the desired result.
	\qed
	
	\begin{lem}\label{widetlide-D_k_0-bound-lemma}
		For $N\in\mathbb{N}$, $a\in\{1,2,\ldots,N\}$, $k_0, d\in\mathbb{N}_0$,  $\boldsymbol{c}=(c_1,\ldots,c_d)\in(\mathbb{N}\setminus\{2\})^d$ with $c_1\geq3$, $\boldsymbol{z}=(z_0,z_1,\ldots,z_d)\in\mathbb{C}^{d+1}$
		with $|z_j|<1, j=0,1,\ldots,d$, let
		\begin{align}\label{Definition of widetlide-D_k_0}
			\widetilde{D}_{k_0}(\boldsymbol{c};\boldsymbol{z})=\sum_{k_0\geq k_1\geq\cdots\geq k_d\geq0}\frac{\left(\frac{2a}{N}\right)_{k_d}}{k_d!(Nk_d+a)}\prod_{i=0}^{d}\frac{(-1)^{k_i\delta_i}(Nk_i+a)^{\delta_i-1}}{(Nk_i+a)^2-z_i^2}W_{k_{i-1},k_i}^{\#}(\{1\}^{c_i-3}),
		\end{align}
		where $k_{-1}=-1$, $\delta_i=\delta(c_i)+\delta(c_{i+1})$ with $c_0=1$ and $c_{d+1}=0$. Then there exist $\varepsilon, p, q, C \in\mathbb{R}$ with $p> 1$ and $\varepsilon, C>0$ such that for all $k_0\in\mathbb{N}_0$,
		$$\left|\widetilde{D}_{k_0}(\boldsymbol{c};\boldsymbol{z})\right|< \frac{C\log^q(Nk_0+N+a)}{k_0^p+\varepsilon}.$$
	\end{lem}
	
	\proof
	For $i=1,2,\ldots,d$, if $c_i\in\{1, 3\}$,  we have
	\begin{align*}
		W_{k_{i-1},k_i}^{\#}(\{1\}^{c_i-3})=2^{\triangle(k_{i-1},k_i)}\leq2,
	\end{align*}
	and if $c_i>3$, we have
	\begin{align*}
		W_{k_{i-1},k_i}^{\#}(\{1\}^{c_i-3})&=\sum_{k_{i-1}\geq l_1\geq\cdots\geq l_{c_i-3}\geq k_i}\frac{2^{\triangle(k_{i-1},l_1)+\triangle(l_1,l_2)+\cdots+\triangle(l_{c_i-3},k_i)}}{(Nl_1+a)\cdots(Nl_{c_i-3}+a)}\\
		&<C_{1}\log^{c_i-3}(Nk_{i-1}+N+a),
	\end{align*}
	where $C_{1}>2$ is some positive constant. Then we can express the above in a unified form
	\begin{align*}
		W_{k_{i-1},k_i}^{\#}(\{1\}^{c_i-3})
		<C_2{\log}^{c_i-1}(Nk_{0}+N+a),
	\end{align*}
	where $C_2>2(\log2)^{-2}$ is a positive constant. Notice that $\delta_0=1$, $\delta_d\in\{2, 3\}$, $\delta_i\in\{0,1,2\}$ for $i=1,2,\ldots,d-1$, and $$\frac{\left(\frac{2a}{N}\right)_{k_d}}{k_d!(Nk_d+a)}<\frac{k_d+1}{Nk_d+a}\leq1.$$
	So we deduce that
	\begin{align*}
		&\left|\widetilde{D}_{k_0}(\boldsymbol{c};\boldsymbol{z})\right|\\
		=&\frac{2}{|(Nk_0+a)^2-z_0^2|}\left|\sum_{k_0\geq k_1\geq\cdots\geq k_d\geq0}
        \frac{\left(\frac{2a}{N}\right)_{k_d}}{k_d!(Nk_d+a)}\prod_{i=1}^{d}\frac{(-1)^{k_i\delta_i}(Nk_i+a)^{\delta_i-1}}{(Nk_i+a)^2-z_i^2}W_{k_{i-1},k_i}^{\#}(\{1\}^{c_i-3})\right|\\
		<&C_{3}\cdot\frac{\log^{c_1+\cdots+c_d-d}(Nk_0+N+a)}{(Nk_0+a)^2-|z_0|^2}\sum_{k_0\geq k_1\geq\cdots\geq k_d\geq0}\prod_{i=1}^{d}\frac{Nk_i+a}{(Nk_i+a)^2-|z_i|^2},
	\end{align*}
	where $C_3$ is a positive constant. Let $|z_{\max}|=\max\{|z_1|,\ldots,|z_d|\}$, then
	\begin{align*}
		\left|\widetilde{D}_{k_0}(\boldsymbol{c};\boldsymbol{z})\right|
		&<C_{4}\cdot\frac{\log^{c_1+\cdots+c_d-d}(Nk_0+N+a)}{(Nk_0+a)^2-|z_0|^2}\left(\sum_{k=1}^{k_0}\frac{1}{Nk+a-|z_{\max}|^2}\right)^d\\
		&<C_{5}\cdot\frac{\log^{c_1+\cdots+c_d}(Nk_0+N+a)}{(Nk_0+a)^2-|z_{0}|^2}\\
		&<C_{6}\cdot\frac{\log^{c_1+\cdots+c_d}(Nk_0+N+a)}{k_0^2+a^2-|z_{0}|^2},
	\end{align*}
	where $C_{4}, C_{5}$ and $C_{6}$ are positive constants. Hence by setting $\varepsilon=a^2-|z_{0}|^2$, $q=c_1+\cdots+c_d$, $p=2$ and $C=C_6$, we conclude the result.
	\qed
	
	\subsection{Main theorems}\label{SubSec:ProfThms3-4,5}
	
		The following theorem generalizes the results of \cite[Theorem 1.4]{PP2019} and \cite[Theorem 3.1]{Li-Yan}.
	
	\begin{thm}\label{D-main-thm}
		Let $a, N$ be integers with $1\leq a \leq N$. For any $d\in\mathbb{N}_0$, $\boldsymbol{c}=(c_1,\ldots,c_d)\in(\mathbb{N}\setminus\{2\})^d$ with $c_1\geq3$, and $\boldsymbol{z}=(z_0,z_1,\ldots,z_d)\in\mathbb{C}^{d+1}$ with $|z_j|<1, j=0,1,\ldots,d$, we have
		\begin{align}\label{D-main-formula}
			&D(\boldsymbol{c};\boldsymbol{z})\nonumber\\
			=&\frac{N}{B\left(\frac{a}{N},\frac{a}{N}\right)}\sum_{k_0\geq k_1\geq\cdots\geq k_d\geq0}\frac{\left(\frac{2a}{N}\right)_{k_d}}{k_d!(Nk_d+a)}
			\prod_{i=0}^{d}\frac{(-1)^{k_i\delta_i}(Nk_i+a)^{\delta_i-1}}{(Nk_i+a)^2-z_i^2}W_{k_{i-1},k_i}^{\#}(\{1\}^{c_i-3}),
		\end{align}
		where $k_{-1}=-1$, $\delta_i=\delta(c_i)+\delta(c_{i+1})$ with $c_0=1$ and $c_{d+1}=0$, $(a)_n$ is the Pochhammer symbol, $B(a,b)$ is the Beta function.
	\end{thm}
	
	The following theorem provides explicit expressions of multiple $t$-star values of level $N$ with an arbitrary number of blocks of twos.
	
	\begin{thm}\label{t^star-levelN-theorem}
		Let $a, N$ be integers with $1\leq a \leq N$. For any $d,b_0,b_1,\ldots,b_d\in\mathbb{N}_0$, $c_1,\ldots,c_d\in\mathbb{N}\setminus\{2\}$ with $c_1\geq3$ if $b_0=0$ and $d\geq1$, we have
		\begin{align}\label{t^star-levelN-formula}
			&t_{N,a}^{\star}(\{2\}^{b_0},c_1,\{2\}^{b_1},\ldots,c_d,\{2\}^{b_d})\nonumber\\
			=&\frac{N}{B\left(\frac{a}{N},\frac{a}{N}\right)}\sum_{k_0\geq k_1\geq\cdots\geq k_d\geq0}\frac{\left(\frac{2a}{N}\right)_{k_d}}{k_d!(Nk_d+a)}
			\prod_{i=0}^{d}\frac{(-1)^{k_i\delta_i}}{(Nk_i+a)^{2b_i+3-\delta_i}}W_{k_{i-1},k_i}^{\#}(\{1\}^{c_i-3}),
		\end{align}
		where $k_{-1}=-1$, $\delta_i=\delta(c_i)+\delta(c_{i+1})$ with $c_0=1$ and $c_{d+1}=0$, $(a)_n$ is the Pochhammer symbol, $B(a,b)$ is the Beta function.
	\end{thm}

	\noindent \textbf{Proof of Theorem \ref{D-main-thm}}
	By Lemma \ref{limit-D_n=D-lemma} and Theorem \ref{D_n-main-thm}, we obtain that
	$$D(\boldsymbol{c};\boldsymbol{z})=\lim_{n\rightarrow\infty}D_n(\boldsymbol{c};\boldsymbol{z})
	=\lim_{n\rightarrow\infty}\sum_{k_0=0}^{n}\frac{N\left[\left(\frac{a}{N}\right)_{n+1}\right]^2}{(n-k_0)!\left(\frac{2a}{N}\right)_{n+k_0+1}}\widetilde{D}_{k_0}(\boldsymbol{c};\boldsymbol{z}),$$
	where $\widetilde{D}_{k_0}(\boldsymbol{c};\boldsymbol{z})$ is defined in \eqref{Definition of widetlide-D_k_0}. Applying Lemmas \ref{limit-M_k-lemma} and \ref{widetlide-D_k_0-bound-lemma}, we have
	$$D(\boldsymbol{c};\boldsymbol{z})
	=\frac{N}{B\left(\frac{a}{N},\frac{a}{N}\right)}\sum_{k_0=0}^{\infty}\widetilde{D}_{k_0}(\boldsymbol{c};\boldsymbol{z}).$$
	This completes the proof.
	\qed
	
	\noindent \textbf{Proof of Theorem \ref{t^star-levelN-theorem}}
	Set $\boldsymbol{b}=(b_0,b_1,$\ldots$,b_d)$ and
	\begin{align*}
		\widetilde{D}_{k_0}(\boldsymbol{c}; \boldsymbol{b})=\sum_{k_0\geq k_1\geq\cdots\geq k_d\geq0}\frac{\left(\frac{2a}{N}\right)_{k_d}}{k_d!(Nk_d+a)}\prod_{i=0}^{d}\frac{(-1)^{k_i\delta_i}}{(Nk_i+a)^{2b_i-\delta_i+3}}W_{k_{i-1},k_i}^{\#}(\{1\}^{c_i-3}).
	\end{align*}
	According to the proof of Lemma \ref{widetlide-D_k_0-bound-lemma}, it is easy to deduce that there exist $\varepsilon, p, q, C \in\mathbb{R}$ with $p> 1$ and $\varepsilon, C>0$ such that for all $k_0\in\mathbb{N}_0$,
	\begin{align}\label{widetilde-D-k_0}
		\left|\widetilde{D}_{k_0}(\boldsymbol{c}; \boldsymbol{b})\right|< \frac{C\log^q(Nk_0+N+a)}{k_0^p+\varepsilon}.
	\end{align}
	By taking the limit $n\rightarrow\infty$ in Corollary \ref{_nt_N,a^star-corollary} and using \eqref{widetilde-D-k_0} and Lemma \ref{limit-M_k-lemma}, we conclude the desired result.
	\qed


	\section{Evaluations for multiple $t$-star values of level $N$}\label{Sec: Evaluations for MtSVs of level $N$}
	
	In this section, we provide some evaluations of multiple $t$-star values of level $N$. For $\boldsymbol{s}=(s_1,\ldots,s_u)\in(\mathbb{Z}\backslash\{0\})^u$ with $s_1\neq 1$, we define the alternating multiple $t$-value of level $N$ containing factorial and pochhammer symbols by
	\begin{align*}
		\widetilde{t}_{N,a}(\boldsymbol{s})=\widetilde{t}_{N,a}(s_1,\ldots,s_u)
		=\sum_{k_1>\cdots>k_u\geq0}\frac{{\sign(s_1)}^{k_1}\cdots{\sign(s_u)}^{k_u}}{{(Nk_1+a)}^{|s_1|}\cdots {(Nk_u+a)}^{|s_u|}}\cdot\frac{\left(\frac{2a}{N}\right)_{k_u}}{k_u!},
	\end{align*}
	where $\sign(s)$ is $1$ if $s>0$ and $-1$ if $s<0$. Let $r$ be a variable. The interpolated alternating multiple $t$-values of level $N$ are defined by
	\begin{align*}
		{\widetilde{t}}_{N,a}^r(\boldsymbol{s})=\widetilde{t}_{N,a}^r(s_1,\ldots,s_u)=\sum_{\circ=\lq\lq+"\text{or}\lq\lq,"}r^{\#\text{pluses}}\widetilde{t}_{N,a}(s_1\circ s_2\circ\cdots\circ s_u),
	\end{align*}
	where
	\[
	s+l=
	\begin{cases}
		|s|+|l| &\text{if\;} sl>0,\\
		-(|s|+|l|) &\text{if\;} sl<0.
	\end{cases}
	\]
As usual, we may use the bar notation. For example, for positive integers $s_1,s_2$, we denote $\widetilde{t}_{N,a}(-s_1,s_2)$ by $\widetilde{t}_{N,a}(\overline{s_1},s_2)$ and denote  $\widetilde{t}_{N,a}^r(-s_1,s_2)$ by $\widetilde{t}_{N,a}^r(\overline{s_1},s_2)$.

	\begin{rem}
		Setting $N=a=1$, we get the interpolated alternating multiple zeta value,
		\begin{align*}
			\zeta^r(\boldsymbol{s})=\zeta^r(s_1,\ldots,s_u)=\sum_{\circ=\lq\lq+"\text{or}\lq\lq,"}r^{\#\text{pluses}}\zeta(s_1\circ s_2\circ\cdots\circ s_u),
		\end{align*}
		which is an alternate version of interpolated multiple zeta value as defined by S. Yamamoto \cite{Yamamoto}. Setting $N=2, a=1$, we obtain the interpolated alternating multiple $t$-values \cite{Li-Yan},
		\begin{align*}
			t^r(\boldsymbol{s})=t^r(s_1,\ldots,s_u)=\sum_{\circ=\lq\lq+"\text{or}\lq\lq,"}r^{\#\text{pluses}}t(s_1\circ s_2\circ\cdots\circ s_u).
		\end{align*}
	\end{rem}
	
	Therefore, Theorem \ref{t^star-levelN-theorem} can be rewritten as the following theorem.
	
		\begin{thm}\label{t^star-levelN-theorem-4.2}
		Let $a, N$ be integers with $1\leq a \leq N$. For any $d,b_0,b_1,\ldots,b_d\in\mathbb{N}_0$, $c_1,\ldots,c_d\in\mathbb{N}\setminus\{2\}$ with $c_1\geq3$ if $b_0=0$ and $d\geq1$, we have
		\begin{align}\label{t^star-levelN-formula-4.2}
			&t_{N,a}^{\star}(\{2\}^{b_0},c_1,\{2\}^{b_1},\ldots,c_d,\{2\}^{b_d})
			=\frac{N\cdot2^{|\boldsymbol{c}|+1-2d+2\delta(\boldsymbol{c})}}{B\left(\frac{a}{N},\frac{a}{N}\right)}\nonumber\\
			&\times{\widetilde{t}}_{N,a}^{1/2}((-1)^{\delta(c_1)-1}(2b_0+2-\delta(c_1)),\{1\}^{c_1-3},(-1)^{\delta(c_1)+\delta(c_2)}(2b_1+3-\delta(c_1)-\delta(c_2)),\nonumber\\
			&\{1\}^{c_1-3},\ldots,(-1)^{\delta(c_{d-1})+\delta(c_{d})}(2b_{d-1}+3-\delta(c_{d-1})-\delta(c_d)),\{1\}^{c_d-3},\nonumber\\
			&(-1)^{\delta(c_{d})}(2b_{d}+2-\delta(c_d))),
		\end{align}
		where $B(a,b)$ is the Beta function, $c_0=1$, $|\boldsymbol{c}|=c_1+\cdots+c_d$ and $\delta(
		\boldsymbol{c})=\delta(c_1)+\cdots\delta(c_d)$ with
		\[
		\delta(c)=
		\begin{cases}
			1 &\text{if\;} c=1,\\
			0 &\text{if\;} c\geq3.
		\end{cases}
		\]
	\end{thm}
	
Taking $N=a=1$, we get formulas for arbitrary multiple zeta-star values $\zeta^{\star}(\{2\}^{b_0},$\\$c_1, \{2\}^{b_1},\ldots,c_d,\{2\}^{b_d})$ in terms of interpolated alternating multiple zeta values\cite[Theorem 1.4]{PPZ2016}. Let $N=2, a=1$, the main result of Li-Yan \cite[Theorem 3.3]{Li-Yan} can be restated as the following cocollary.

\begin{cor}\label{t^star-theorem-4.4}
			For any $d,b_0,b_1,\ldots,b_d\in\mathbb{N}_0$, $c_1,\ldots,c_d\in\mathbb{N}\setminus\{2\}$ with $c_1\geq3$ if $b_0=0$ and $d\geq1$, we have
			\begin{align}\label{t^star-formula-4.4}
				&t^{\star}(\{2\}^{b_0},c_1,\{2\}^{b_1},\ldots,c_d,\{2\}^{b_d})
				=-\frac{2^{|\boldsymbol{c}|+2-2d+2\delta(\boldsymbol{c})}}{\pi}\nonumber\\
				&\times t^{1/2}((-1)^{\delta(c_1)-1}(2b_0+2-\delta(c_1)),\{1\}^{c_1-3},(-1)^{\delta(c_1)+\delta(c_2)}(2b_1+3-\delta(c_1)-\delta(c_2)),\nonumber\\
				&\{1\}^{c_1-3},\ldots,(-1)^{\delta(c_{d-1})+\delta(c_{d})}(2b_{d-1}+3-\delta(c_{d-1})-\delta(c_d)),\{1\}^{c_d-3},\nonumber\\
				&(-1)^{\delta(c_{d})}(2b_{d}+2-\delta(c_d))),
			\end{align}
			with the same notation as in Theorem \ref{t^star-levelN-theorem-4.2}.
		\end{cor}

	Setting $d=0$ in Theorem \ref{t^star-levelN-theorem-4.2}, we obtain the following formula.
	
	\begin{cor}
		For any $b\in\mathbb{N}_0$, we have
		\begin{align*}
			t_{N,a}^{\star}(\{2\}^b)=\frac{2N}{B\left(\frac{a}{N},\frac{a}{N}\right)}	\widetilde{t}_{N,a}(\overline{2b+1}).
		\end{align*}
	\end{cor}
	
	Setting $c_1=\cdots=c_{d}=1$,  $c_1=\cdots=c_{d}=3$ and $c_1, \ldots, c_{d}\geq3$ respectively in Theorem \ref{t^star-levelN-theorem-4.2}, we have the following three corollaries.
	\begin{cor}
	For any $d, b_0 \in\mathbb{N}$, $b_1, \ldots, b_d \in\mathbb{N}_0$, we have
	\begin{align*}
		t_{N,a}^{\star}(\{2\}^{b_0}, 1, \{2\}^{b_1}, \ldots, 1, \{2\}^{b_d})=\frac{2^{d+1}N}{B\left(\frac{a}{N},\frac{a}{N}\right)}\widetilde{t}_{N,a}^{1/2}(2b_0+1, 2b_1+1, \ldots, 2b_{d-1}+1, \overline{2b_d+1}).
	\end{align*}
   \end{cor}

	\begin{cor}
	For any $d, b_0 \in\mathbb{N}$, $b_1, \ldots, b_d \in\mathbb{N}_0$, we have
    \begin{align*}
	t_{N,a}^{\star}(\{2\}^{b_0}, 3, \{2\}^{b_1}, \ldots, 3, \{2\}^{b_d})=\frac{2^{d+1}N}{B\left(\frac{a}{N},\frac{a}{N}\right)}\widetilde{t}_{N,a}^{1/2}(\overline{2b_0+2}, 2b_1+3, \ldots, 2b_{d-1}+3, 2b_d+2).
    \end{align*}
	\end{cor}
	
	\begin{cor}
		For any $d, b_0, b_1, \ldots, b_d \in\mathbb{N}_0$, $c_1, \ldots, c_d\in\mathbb{N}$ with $c_1, \ldots, c_d\geq3$, we have
		\begin{align*}
			&t_{N,a}^{\star}(\{2\}^{b_0}, c_1, \{2\}^{b_1}, \ldots, c_d, \{2\}^{b_d})\\
			&=\frac{2^{c1+\cdots+c_d+1-2d}N}{B\left(\frac{a}{N},\frac{a}{N}\right)}
			\widetilde{t}_{N,a}^{1/2}(\overline{2b_0+2},\{1\}^{c_1-3},2b_1+3, \{1\}^{c_2-3}, \ldots, 2b_{d-1}+3, \{1\}^{c_d-3},2b_d+2).
		\end{align*}
	\end{cor}
	
	Setting $(c_1, c_2, \ldots, c_{2d})=(\{3, 1\}^d)$ and $(c_1, c_2, \ldots, c_{2d})=(\{1, 3\}^d)$ in Theorem \ref{t^star-levelN-theorem-4.2}, we deduce the following results.
	
	\begin{cor}
		For any $d, b_0, \ldots, b_{2d}\in\mathbb{N}_0$, we have
		\begin{align*}
			&t_{N,a}^{\star}(\{2\}^{b_0}, 3, \{2\}^{b_1}, 1, \{2\}^{b_2}, 3, \{2\}^{b_3}, 1 \ldots, \{2\}^{b_{2d-2}}, 3, \{2\}^{b_{2d-1}}, 1,  \{2\}^{b_d})\\
			&=\frac{2^{2d+1}N}{B\left(\frac{a}{N},\frac{a}{N}\right)}\widetilde{t}_{N,a}^{1/2}(\overline{2b_0+2}, \overline{2b_1+2}, \ldots, \overline{2b_{2d-1}+2}, \overline{2b_{2d}+1}).
		\end{align*}
	\end{cor}

	\begin{cor}
	For any $d, b_0\in\mathbb{N}$, $b_1, \ldots, b_{2d}\in\mathbb{N}_0$, we have
	\begin{align*}
		&t_{N,a}^{\star}(\{2\}^{b_0}, 1, \{2\}^{b_1}, 3, \{2\}^{b_2}, 1, \{2\}^{b_3}, 3 \ldots, \{2\}^{b_{2d-2}}, 1, \{2\}^{b_{2d-1}}, 3,  \{2\}^{b_d})\\
		&=\frac{2^{2d+1}N}{B\left(\frac{a}{N},\frac{a}{N}\right)}\widetilde{t}_{N,a}^{1/2}(2b_0+1, \overline{2b_1+2}, \ldots, \overline{2b_{2d-1}+2}, 2b_{2d}+2).
	\end{align*}
    \end{cor}
	
	\vskip10pt
	
	\noindent{\bf Acknowledgements.} This work is supported by the Fundamental Research Funds for the Central Universities (grant number 22120210552).

	\medskip



\begin{thebibliography}{99}	
		
		
		
		\bibitem{BBBL1998}
		J. M. Borwein, D. M. Bradley, D. J. Broadhurst and P. Lison\v{e}k, Combinatorial aspects of multiple zeta values, \emph{Electron. J. Comb.} \textbf{5} (1998), R38.
		
		\bibitem{BBBL2001}
		J. M. Borwein, D. M. Bradley, D. J. Broadhurst and P. Lison\v{e}k, Special values of multiple polylogarithms, \emph{Trans. Amer. Math. Soc.} \textbf{353} (2001), no. 3, pp. 907--941.
		
		\bibitem{SCharlton}
		S. Charlton, On motivic multiple $t$ values, Saha's basis conjecture, and generators of alternating MZV's, preprint, arXiv: 2112.14613.
		
		\bibitem{Chung2019}
		C. Chung, On the sum relation of multiple Hurwitz zeta functions, \emph{Quaest. Math.} \textbf{42} (2019), no. 3, pp. 297--305.
		
		\bibitem{Diaconis-Freedman}
		P. Diaconis and D. Freedman, An elementary proof of Stirling’s formula, \emph{Amer. Math. Monthly} \textbf{93} (1986), pp. 123--125.
		
		\bibitem{Euler}	
		L. Euler, \emph{Institutiones calculi differentialis}, G. Kowalewski (ed.), Opera Omnia Ser. 1; opera mat., \textbf{10}, Teubner (1980).
		
		\bibitem{Feller}
		W. Feller, A direct proof of Stirling’s formula, \emph{Amer. Math. Monthly} \textbf{74} (1967), pp.1223--1225.
		
		\bibitem{PP2017}
		Kh. Hessami Pilehrood and T. Hessami Pilehrood, An alternative proof of a theorem of Zagier, \emph{J. Math. Anal. Appl.} \textbf{449} (2017), no. 1, pp. 168--175.
		
		\bibitem{PP2019}
		Kh. Hessami Pilehrood and T. Hessami Pilehrood, Generating functions for multiple zeta star values, \emph{J. Th\'{e}or. Nombres Bordeaux} \textbf{31} (2019), no. 2, pp. 343--360.
		
		\bibitem{PP2022}
		Kh. Hessami Pilehrood and T. Hessami Pilehrood, Multiple zeta star values on $3-2-1$ indices, \emph{Ramanujan J} \textbf{60} (2022), pp. 259--285.
		
		\bibitem{PPT2014}
		Kh. Hessami Pilehrood, T. Hessami Pilehrood and R. Tauraso, New properties of multiple harmonic sums modulo $p$ and $p$-analogues of Leshchiner's series, \emph{Trans. Amer. Math. Soc.} \textbf{366} (2014), no. 6, pp. 3131--3159.
		
		\bibitem{PPZ2016}
		Kh. Hessami Pilehrood, T. Hessami Pilehrood and J. Zhao, On $q$-analogs of some families of multiple harmonic sums and multiple zeta star value identities, \emph{Commun. Number Theory Phys.} \textbf{10} (2016), no. 4, pp. 805--832.
		
		\bibitem{Hoffman1992}
		M. E. Hoffman, Multiple harmonic series, \emph{Pacific J. Math.} \textbf{152} (1992), no. 2, pp. 275--290.
		
		\bibitem{Hoffman2019}
		M. E. Hoffman, An odd variant of multiple zeta values, \emph{Comm. Number Theory Phys.} \textbf{13} (2019), no. 3, pp. 529--567.
		
		\bibitem{Lee-Peng2018}		
		T. Lee-Peng, Alternating double Euler sums, hypergeometric identities and a theorem of Zagier, \emph{J. Math. Anal. Appl.} \textbf{462} (2018), no. 1, pp. 777--800.	
		
		\bibitem{Li2013}
		Z. Li, Another proof of Zagier's evaluation formula of the multiple zeta values $\zeta (2, \dots , 2, 3, 2, \dots , 2)$, \emph{Math. Res. Lett.} \textbf{20} (2013), no. 5, pp. 947--950.
		
		\bibitem{Li-Lupu-Orr}
		L. Lai, C. Lupu and D. Orr, Elementary proofs of Zagier's formula for multiple zeta values and its odd variant, preprint, arXiv: 2201.09262.
		
		\bibitem{Li-Wang}
		Z. Li and Z. Wang, Relations of multiple $t$-values of general level, preprint, arXiv: 2210.16854.
		
		\bibitem{Li-Yan}
		Z. Li and L. Yan, Generating functions of multiple $t$-star values, \emph{Quaest. Math.}, accepted.
		
		\bibitem{Linebarger-Zhao}	
		E. Linebarger and J. Zhao, A family of multiple harmonic sum and multiple zeta star value identities, \emph{Mathematika} \textbf{61} (2015), no. 1, pp. 63--71.
		
		\bibitem{SMuneta2008}
		S. Muneta, On some explicit evaluations of multiple zeta-star values, \emph{J. Number Theory} \textbf{128} (2008), no. 9, pp. 2538--2548.
		
		\bibitem{TMurakami}
		T. Murakami, On Hoffman's $t$-values of maximal height and generators of multiple zeta values, \emph{Math. Ann.} \textbf{382} (2022), no. 1-2, pp. 421--458.
		
		\bibitem{Ohno-Zudilin}
		Y. Ohno and W. Zudilin, Zeta stars, \emph{Commun. Number Theory Phys.} \textbf{2} (2008), no. 2, pp. 325--347.
		
		\bibitem{Quan2020}
		J. Quan, Alternating double $t$-values and $T$-values, \emph{Adv. Differ. Equ.} \textbf{2020} 450 (2020).
		
		\bibitem{Wallis1656}
		J. Wallis, \emph{Arithmetica Infinitorum}, Oxford, England, 1656.
		
		\bibitem{Xu-Yan}
		C. Xu and L. Yan, Parametric Euler $T$-sums of odd harmonic numbers, preprint, arXiv: 2203.13996.
		
		\bibitem{Yamamoto}	
		S. Yamamoto, Interpolation of multiple zeta and zeta-star values, \emph{J. Algebra} \textbf{385} (2013), pp. 102--114.		
				
		\bibitem{Yuan-Zhao}
		H. Yuan and J. Zhao, Double shuffle relations of double zeta values and double Eisenstein series of level $N$, \emph{J. London Math. Soc.} \textbf{92} (2) (2015), pp. 520--546.
		
		\bibitem{DZagier1994}
		D. Zagier, \emph{Values of zeta functions and their applications}, First European Congress of Mathematics, vol. II (Paris, 1992), pp. 497--512, Progr. Math., vol. 120, Birkh\"{a}user, Basel, 1994.
		
		\bibitem{DZagier2012}
		D. Zagier, Evaluation of the multiple zeta values $\zeta(2, \ldots, 2, 3, 2, \ldots , 2)$, \emph{Ann. Math.} \textbf{175} (2012), no. 2, pp. 977--1000.
		
		\bibitem{Zhao2016}
		J. Zhao, Identity families of multiple harmonic sums and multiple zeta star values, \emph{J. Math. Soc. Japan} \textbf{68} (2016), no. 4, pp. 1669--1694.
		
	\end{thebibliography}
\end{document}